\documentclass[A4,draftcls,onecolumn,10pt]{IEEEtran}
\usepackage{url}
\usepackage[dvips]{graphicx}
\usepackage{amssymb,amsfonts,amsmath,amsthm,amscd,dsfont,mathrsfs}
\usepackage{graphicx,float,psfrag,epsfig,color}
\usepackage{subfig}
\usepackage{cite}
\usepackage{setspace}
\newcommand*{\Scale}[2][4]{\scalebox{#1}{$#2$}}%

\newcommand{\R}{\mathbb{R}} 
\newcommand{\N}{\mathbb{N}}

\newcommand{\pder}[1]{{\frac{\partial}{\partial #1}}}
\newcommand{\sdet}{{\mathsf{det}}}
\newcommand{\strace}{{\mathsf{tr}}}
\newcommand{\sadj}{{\mathsf{adj}}}
\newcommand{\semb}{{\mathsf{embb}}}
\newcommand{\inp}[2]{{\langle #1, #2 \rangle}}

\newcommand{\part}{\mathcal{P}}

\theoremstyle{definition}
\theoremstyle{example}
\theoremstyle{plain}
\newtheorem{thm}{Theorem}
\theoremstyle{plain}
\theoremstyle{plain}
\newtheorem{lemma}{Lemma}
\theoremstyle{plain}
\theoremstyle{plain}
\theoremstyle{remark}
\theoremstyle{discussion}
\theoremstyle{plain}

\begin{document}

\title{A New Identity for the Least-square Solution of Overdetermined Set of Linear Equations}

\author{\IEEEauthorblockN{Saeid Haghighatshoar\IEEEauthorrefmark{1},
Mohammad J. Taghizadeh\IEEEauthorrefmark{2}, and
Afsaneh Asaei\IEEEauthorrefmark{3}}\\
\IEEEauthorrefmark{1}Technische Universit\"at, Berlin, Germany\\ 
\IEEEauthorrefmark{2}Huawei European Research Centre, Munich, Germany\\ 
\IEEEauthorblockA{\IEEEauthorrefmark{3}Idiap Research Institute, Martigny, Switzerland\\
{\small Emails: saeid.haghighatshoar@tu-berlin.de, mohammad.taghizadeh@huawei.com, afsaneh.asaei@idiap.ch}}}
\maketitle 

\begin{abstract}
In this paper, we prove a new identity for the least-square solution of an over-determined set of linear equation $Ax=b$, where $A$ is an $m\times n$ full-rank matrix, $b$ is a column-vector of dimension $m$, and $m$ (the number of equations) is larger than or equal to  $n$ (the dimension of the unknown vector $x$). Generally, the equations are inconsistent and there is no feasible solution for $x$ unless $b$ belongs to the column-span of $A$. In the least-square approach, a candidate solution is found as the unique $x$ that minimizes the error function $\|Ax-b\|_2$. 

We propose a more general approach that consist in considering all the consistent subset of the equations, finding their solutions, and   taking a weighted average of them to build a candidate solution. In particular, we show that by weighting the solutions with the squared determinant of their coefficient matrix, the resulting candidate solution  coincides with the least square solution.  
\end{abstract}


\section{Introduction}
\subsection{Over-determined Set of Linear Equations}
Let $A$ be an $m\times n$ full-rank matrix and let $b \in \R^m$ be a column vector, and consider the linear equation $Ax=b$, to be solved for the unknown vector $x \in \R^n$. Theory and practice of solving these equations play a major role in essentially every part of mathematics such as linear algebra, operational research, optimization, combinatorics, etc. When $m >n$, we call the equations over-determined and there is a solution if and only if $b$ belongs to the column-span  of $A$~\cite{LS2}. Generally, the equations are inconsistent and we need some kind of criteria to build a candidate solution. 

One approach for finding a solution is the least-square approach~\cite{wikipedia}, where we find a solution by minimizing the quadratic form $\|Ax-b\|_2^2$. The resulting solution is given by $\hat{x}=A^{\#} b$, where $A^{\#}=(A^tA)^{-1} A^t$ denotes the pseudo-inverse of $A$. In estimation theory, $\hat{x}$ can be interpreted as the \textit{best linear unbiased  estimator} (BLUE) of the signal $x$ observed via a linear channel given by the matrix $A$ and contaminated with an i.i.d. Gaussian noise~\cite{LS}. 
Note that in this case, if $b$ is in the column-span of $A$, the resulting estimation error is zero.

Another approach for building a candidate solution is by some kind of averaging all the possible sub-solutions. To explain this more precisely, we first need to introduce some notations. For $k \in \N$, we define $[k]=\{1,2,\dots,k\}$ to be the set of all integers from $1$ up to $k$. 
We denote by $\part$  the set of all subsets of $[m]$ of size $n$, i.e., $\part=\{p \subset [m]: |p|=n\}$, where $|p|$ denotes the size of the subset $p$. For a $p\subset \part$, we define $A_p$ to be the $n\times n$ matrix obtained by selecting the rows of the matrix $A$ belonging to $p$ by keeping their order as in $A$. 

Suppose $p \in \part$ is such that $\sdet(A_p) \neq 0$. By restricting the equations to $A_p$, we can obtain a sub-solution $x_p=A_p^{-1} b_p$, where $b_p$ is the a sub-vector of $b$ consisting of the components with index in $p$ whose order is the same as in $b$. Taking the weighted average of all possible sub-solutions with a weighting $\omega_p\geq 0$, $p\in \part$, we can build a candidate solution as follows
\begin{align}
s^\omega=\frac{\sum_{p \in \part} \omega_p x_p}{\sum_{p \in \part} \omega_p}.
\end{align}
As the matrix $A$ is full-rank, there is at least one $p \in \part$ with a nonzero $\sdet(A_p)$, thus $s^\omega$ is well-defined.
By changing the associated weighting $\omega_p$, we obtain a variety of candidate solutions for the over-determined equation $Ax=b$. 

Let us consider the weighting function $\omega_p=\sdet(A_p)^2$, which is equal to the squared determinant of the sub-matrix $A_p$, and let us define the resulting solution by 
\begin{align}\label{w_a_s}
\hat{x}_{\mathsf{LS}}=\frac{\sum_{p \in \part} \sdet(A_p)^2 A_p^{-1} b_p}{\sum_{p \in \part} \sdet(A_p)^2}.
\end{align}
If for a specific $p \in \part$, $\sdet(A_p)=0$ then $A_p^{-1}$ does not exist but, with some abuse of notation, this term does not play a role because its corresponding weight $\sdet(A_p)^2$ is equal to $0$.

\subsection{Our Contribution}
We prove that with the weighting $\omega_p=\sdet(A_p)^2$, the resulting solution $\hat{x}_{\mathsf{LS}}$ in Eq.~\eqref{w_a_s} coincides with the least-square solution given by $A^{\#} b=(A^t A)^{-1} A^t b$. More importantly, this holds for every full-rank matrix $A$ and for an arbitrary vector $b$. We have summarized this in the following theorem.
\begin{thm}\label{main_thm}
Suppose $A$ is a given $m\times n$ full-rank matrix  with $m\geq n$ and assume that $b \in \R^m$ is an arbitrary vector. Let $\hat{x}_{\mathsf{LS}}$ be the weighted average solution given by Eq.~\eqref{w_a_s}. Then $\hat{x}_{\mathsf{LS}}=(A^t A)^{-1} A^t b$, i.e., $\hat{x}_{\mathsf{LS}}$ coincides with the least square solution.
\end{thm}

\subsection{Notation and Auxiliary Results}
In this section, we first introduce the required notations for the rest of the paper and prove some auxiliary results that we need to prove Theorem \ref{main_thm}. Let $B$ be an arbitrary $n \times n$ matrix and let $p \subset [m]$ of size $|p|=n$. We denote by $\semb(B,p,m)$ the embedding of columns of $B$ inside an $n\times m$ matrix. More precisely, assume that the components of $p$ are sorted with $p_1 <p_2<\dots<p_n$. Then $\semb(B,p,m)$ is an $n \times m$ matrix whose $p_i$-th column, $i\in [n]$, is equal to the $i$-th column of $B$, and all the other $m-n$ columns are set to zero. 
 
Let $r,c \in \N$ be arbitrary numbers. We define the linear space of all $r\times c$ real-valued matrices by  $M_\R(r,c)$ with the traditional matrix addition and scalar-matrix multiplication. For arbitrary matrices $M,N \in M_\R(r,c)$, we define the following bilinear form $\inp{M}{N}=\strace(M N^t)=\sum_{i,j} M_{ij} N_{ij}$. It is not difficult to see that $\inp{}{}$ defines an inner product on $M_\R(r,c)$. We denote the trace and the determinant of a square matrix $M$ by $\strace(M)$ and $\sdet(M)$ respectively.
 We need the following auxiliary results from linear algebra. We have included all the proofs in Appendix \ref{app}.
\begin{lemma}\label{lem:in_zero}
Let $r,c \in \N$ and let $M \in M_\R(r,c)$. If $\inp{M}{N}=0$ for every $N \in M_\R(r,c)$, then $M=0$.
\end{lemma}

\begin{lemma}\label{lem:der_inv}
Let $M$ be an square invertible matrix whose components depend on a parameter $u$. Then, $\pder{u} M^{-1}=-M^{-1} (\pder{u} M) M^{-1}$.
\end{lemma}

\begin{lemma}\label{lem:der_det}
Let $A$ be an square matrix whose components depend on a parameter $u$. Then,
$\pder{u} \sdet(A)=\sdet(A)\strace(A^{-1} \pder{u} A)$
\end{lemma}

\begin{lemma}\label{lem:sym}
Let $M$ and $S$ be $n\times n$ matrices, where $S$ is symmetric. Then $\strace(S M)=\strace( S M^t)$.
\end{lemma}

\begin{thm}[Cauchy-Binet]\label{thm:CB}
Let $A$ and $B$ be $m \times n$ matrices with $m\geq n$. Then,
\begin{align}
\sdet(A^t B)=\sum_{p \subset [m], |p|=n} \sdet(A_p)\sdet(B_p),
\end{align}
where $|p|$ denotes the number of elements of $p \subset [m]$.
\end{thm}

\section{Proof of the Main Theorem}
In the section, we prove Theorem \ref{main_thm}. Using Eq.~\eqref{w_a_s}, we can write $\hat{x}_{\mathsf{LS}}$ in the following form:
\begin{align}\label{eq:s1_form}
\hat{x}_{\mathsf{LS}}(A,b)&=\frac{\sum_{p \in \part} \sdet(A_p)^2 A_p^{-1} b_p}{\sdet(A^t A)}\\
&=\frac{\sum_{p \in \part} \sdet(A_p)^2 \,\semb(A_p^{-1}, p,m) b}{\sdet(A^t A)},
\end{align}
where in the last term we used the definition of $\semb(A_p^{-1}, p,m)$.
Recall that for $p \in \part$, with elements $p_1 < p_2 < \dots < p_n$, we denote by $\semb(A_p^{-1},p,m)$ an all-zero $n \times m$ matrix except for its $p_i$-th column witch is equal to the $i$-th column of $A_p^{-1}$. Now, we need to prove that for any $b \in \R^m$ and for any $m \times n$ full-rank matrix $A$, the following identity holds 
\begin{align}
(A^t A)^{-1} A^t b= \frac{\sum_{p \in \part} \sdet(A_p)^2 \,\semb(A_p^{-1}, p,m)}{\sdet(A^t A)} b.
\end{align}
As this should be true for every $b \in \R^m$, we need to prove  the following matrix identity:
\begin{align}\label{eq:tbp}
\Scale[0.94]{\sdet(A^t A) \,(A^t A)^{-1} A^t = \sum_{p \in \part} \sdet(A_p)^2 \,\semb(A_p^{-1}, p,m).}
\end{align}
As a first step, it is easy to check that both sides are $n \times m$ matrices, thus the dimensions are compatible. 

In order to prove the identity \eqref{eq:tbp}, let us define the function $f: M_\R(m,n) \to \R$ as follows:
\begin{align}\label{eq:f_def}
f(A)=\sdet(A^t A)-\sum_{p \in \part} \sdet(A_p)^2.
\end{align}
Using the Cauchy-Binet formula as stated in Theorem \ref{thm:CB}, we obtain
\begin{align}\label{eq:CB}
\sdet(A^t A)=\sum_{ p \in \part} \sdet(A_p)\sdet(A_p^t)=\sum_{ p \in \part} \sdet(A_p)^2,
\end{align}
which implies that $f(A)=0$ for every $A \in M_\R(m,n)$. Let $u=A_{ij}$ be a parameter denoting the component of $A$ at row $i$ and column $j$. As $f(A)=0$, we have $\pder{u} f(A)=0$, which implies that
\begin{align*}
\pder{u} \sdet(A^t A) &\stackrel{(a)}{=} \sdet(A^t A) \strace\Big\{(A^tA)^{-1} \pder{u}(A^t A)\Big\} \\
&\stackrel{(b)}{=} \sdet(A^t A)\strace\Big\{(A^tA)^{-1} \big ( (\pder{u}A)^t A + A^t \pder{u} A\big )\Big\}\\
&\stackrel{(c)}{=}\sdet(A^t A)\strace\Big\{(A^tA)^{-1} \big ( A^t \pder{u}A + A^t \pder{u} A\big )\Big\}\\
&= 2\, \sdet(A^t A)\strace\Big\{(A^tA)^{-1} A^t \pder{u}A \Big\}\\
&\stackrel{(d)}{=}2\, \sdet(A^t A)\strace\Big\{(A^tA)^{-1} A^t U_{ij} \Big\}\\
&\stackrel{(e)}{=}2\, \sdet(A^t A)\Big \langle(A^tA)^{-1} A^t,U_{ij}^t\Big \rangle,
\end{align*}
where $(a)$ follows from Lemma \ref{lem:der_det} applied to the matrix $A^t A$, $(b)$ follows from the chain rule applied to $A^t A$, $(c)$ follows from Lemma \ref{lem:sym} applied to the symmetric matrix $(A^tA)^{-1}$ and the matrix $(\pder{u}A)^t A$, $(d)$ results by taking the component-wise derivative of $A$ with respect to $u=A_{ij}$ which we denote by $U_{ij}$, and where $(e)$ results from the definition of the inner product for two matrices. We can simply check that $U_{ij}$ is an $m \times n$ matrix with all-zero components except for $ij$-th component which is equal to $1$. 

Now, taking the derivative of the other term in Eq.~\eqref{eq:f_def} with respect to $u=A_{ij}$, we obtain

\begin{align*}
\pder{u} \sum_{p \in \part} \sdet(A_p)^2&= \sum_{p \in \part} 2\, \sdet(A_p) \pder{u} \sdet(A_p)\\
&\stackrel{(a)}{=}\sum_{p \in \part} 2\, \sdet(A_p) \sdet(A_p)\strace(A_p^{-1} \pder{u} A_p)\\
&\stackrel{(b)}{=}\sum_{p \in \part} 2\, \sdet(A_p)^2 \strace(\semb(A_p^{-1},p,m) \pder{u} A)\\
&\stackrel{(c)}{=}2\,\strace\Big \{ \sum_{p \in \part}  \sdet(A_p)^2 \semb(A_p^{-1},p,m) U_{ij} \Big \}\\
&\stackrel{(d)}{=}2\,\Big \langle \sum_{p \in \part} \sdet(A_p)^2 \semb(A_p^{-1},p,m),U_{ij}^t\Big \rangle,
\end{align*}
where $(a)$ results from Lemma \ref{lem:der_det} applied to the matrix $A_p$. We also have $(b)$ from the definition of the embedding $n$ columns of $A_p^{-1}$ in an $m \times n$ matrix. In particular, notice that as the remaining columns of $\semb(A_p^{-1},p,m)$ are all zero, we can replace $A_p$ by $A$. Finally, $(c)$ results from the linearity of the trace operator $\strace$, and $(d)$ follows from the definition of the inner product. Therefore, we obtain that
\begin{align}\label{eq:fin_result}
0&=\pder{u} f(A)= 2 \Big \langle \ U_{ij}^t,\\
&\sdet(A^t A) (A^tA)^{-1} A^t -\sum_{p \in \part} \sdet(A_p)^2 \semb(A_p^{-1},p,m)\, \Big \rangle.\nonumber
\end{align}
Notice that equality in Eq.~\eqref{eq:fin_result} holds for all matrices $U_{ij}^t$, $i \in [m], j\in [n]$. As, $U_{ij}^t$ form an orthonormal basis for the linear space $M_\R(m, n)$, from Lemma \ref{lem:in_zero}, it immediately results that 
\begin{align*}
\Scale[1]{\sdet(A^t A) (A^tA)^{-1} A^t = \sum_{p \in \part} \sdet(A_p)^2 \semb(A_p^{-1},p,m).}
\end{align*}
From Eq.~\eqref{eq:tbp}, this is exactly what we needed to prove.

\begin{appendices}
\section{Proof of the Auxiliary Results}\label{app}
In this section, we provide the proofs of the auxiliary results.

\vspace{1mm}
{\bf Proof of Lemma \ref{lem:in_zero}:} Let $i\in[r], j\in [c]$ be arbitrary numbers and let $N$ be an all-zero matrix except for the $ij$-th element which is set to $1$. It results that 
$$0=\inp{M}{N}=\sum_{k,\ell} M_{k\ell} N_{k\ell}=M_{ij}=0.$$
As this is true for arbitrary $i$ and $j$, it results that $M=0$.

\vspace{1mm}
{\bf Proof of Lemma \ref{lem:der_inv}:}
Let $I$ be the identity matrix of the same order as $M$. Taking derivative from both sides of the identity $I=M M^{-1}$, and using the chain rule, we obtain that $$0=\pder{u} M M^{-1} + M \pder{u} M^{-1},$$ which implies that $\pder{u} M^{-1}=- M^{-1} (\pder{u} M) M^{-1}$.

\vspace{1mm}
{\bf Proof of Lemma \ref{lem:der_det}:}
Assume that $A$ is a $d \times d$ matrix and let us denote by $A_{ij}$ the component of $A$ in row $i$ and column $j$. We first find $\pder{A_{ij}} \sdet(A)$ and use the chain rule to obtain 
\begin{align}\label{chain_rule}
\pder{u} \sdet(A)=\sum_{i,j \in [d]} \pder{A_{ij}} \sdet(A) \pder{u}A_{ij}.
\end{align}
Notice that in order to compute $\sdet(A)$, we can expand it with respect to the $i$-th row, where we obtain 
\begin{align}\label{summation}
\sdet(A)=\sum_{k \in [d]} (-1)^{i+k} \sdet(\tilde{A}_{ik}),
\end{align}
where $\tilde{A}_{ik}$ is a $(d-1)\times (d-1)$ matrix obtained after removing the $i$-th row and the $k$-th column of the matrix $A$. In particular, it can be immediately checked that the only term in the summation \eqref{summation} that depends on $A_{ij}$  is $(-1)^{i+j} \sdet(\tilde{A}_{ij})$, thus we obtain
\begin{align}
\pder{A_{ij}} \sdet(A)= (-1)^{i+j} \sdet(\tilde{A}_{ij})=\sadj(A)_{ji},
\end{align}
where $\sadj(A)$ denotes the \textit{adjoint} of the matrix $A$. Moreover, from the formula $A^{-1}=\frac{\sadj(A)}{\sdet(A)}$ for the inverse of the matrix $A$, we immediately obtain that 
\begin{align}
\pder{A_{ij}} \sdet(A)= \sdet(A) (A^{-1})_{ji}.
\end{align}
Using the the chain-rule as in Eq.~\eqref{chain_rule}, we have
\begin{align*}
\pder{u} \sdet(A)&= \sdet(A) \sum_{ij} (A^{-1})_{ji} \pder{u}A_{ij}=\sdet(A)\strace(A^{-1} \pder{u} A),
\end{align*}
where $\strace$ denotes the trace operator and where $\pder{u} A$ denotes the component-wise partial derivative of $A$ with respect to $u$.

\vspace{1mm}
{\bf Proof of Lemma \ref{lem:sym}:}
The proof simply follows from the properties of the trace operator:
\begin{align*}
\strace(S M)=\strace((S M)^t)=\strace(M^t S^t)=\strace(M^t S)=\strace(S M^t),
\end{align*}
where we used the symmetry of $S$ and the fact that for arbitrary square matrices $K, L$ of the same dimension, $\strace(K L)=\strace(L K)$. 

\end{appendices}

\end{document}